\documentclass[12pt]{amsart}
\usepackage{amssymb,amsmath,amsthm,mathrsfs,MnSymbol}
\usepackage{graphicx}
\usepackage[table]{xcolor}
\usepackage{multirow,scalerel,arydshln,caption,subcaption,enumitem}   
\usepackage{hyperref}
\usepackage[msc-links]{amsrefs}
\hypersetup{
  colorlinks   = true,	
  urlcolor     = blue,
  linkcolor    = blue,
  citecolor   = red
}
\calclayout
\allowdisplaybreaks
\newtheorem{theorem}{Theorem}

\newtheorem{corollary}[theorem]{Corollary}

\newtheorem{proposition}[theorem]{Proposition}
\numberwithin{theorem}{section} \numberwithin{equation}{section}
\newcommand{\beq}{\begin{small} \begin{equation}}
\newcommand{\eeq}{\end{equation} \end{small}}
\newcommand{\beqn}{\begin{small} \begin{equation*}}
\newcommand{\eeqn}{\end{equation*} \end{small}}

\newcommand*{\MyScale}{0.7}
\begin{document}
\title[F-theory/Heterotic String Duality with two Wilson lines]{The duality between F-theory and the Heterotic String in $D=8$ with two Wilson lines}
\author{Adrian Clingher}
\address{Dept.\!~of Mathematics \& Statistics, University of Missouri - St. Louis, St. Louis, MO 63121}
\email{clinghera@umsl.edu}
\thanks{A.C. acknowledges support from a UMSL Mid-Career Research Grant.}
\author{Thomas Hill}
\thanks{T.H. acknowledges the support from the Office of Graduate Studies at Utah State University.}
\author{Andreas Malmendier}
\address{Dept.\!~of Mathematics, University of Connecticut, Storrs, Connecticut 06269\\
\hspace*{0.3cm} Dept.\!~of Mathematics \& Statistics, Utah State University, Logan, UT 84322}
\email{andreas.malmendier@usu.edu}
\thanks{A.M. acknowledges support from the Simons Foundation through grant no.~202367.}
\begin{abstract}
We construct non-geometric string compactifications by using the F-theory dual of the heterotic string compactified on a two-torus with two Wilson line parameters, together with a close connection between modular forms and the equations for certain K3 surfaces of Picard rank $16$. We construct explicit Weierstrass models for all inequivalent Jacobian elliptic fibrations supported on this family of K3 surfaces and express their  parameters in terms of modular forms generalizing Siegel modular forms. In this way, we find a complete list of all dual non-geometric compactifications obtained by the partial higgsing of the heterotic string gauge algebra using two Wilson line parameters.
\end{abstract}
\maketitle
\section{Introduction}
\label{sec:physics}
In a standard compactification of the type IIB string theory, the axio-dilaton field $\tau$ is constant and no D7-branes are present.  Vafa's idea in proposing F-theory \cite{MR1403744} was to simultaneously allow a variable axio-dilaton field $\tau$ and D7-brane sources, defining at a new class of models in which the string coupling is never weak. These compactifications of the type IIB string in which the axio-dilaton field varies over a base are referred to as \emph{F-theory models}. They depend on the following key ingredients:  an $\operatorname{SL}_2(\mathbb{Z})$ symmetry of the physical theory, a complex scalar field $\tau$ with positive imaginary part on which $\operatorname{SL}_2(\mathbb{Z})$ acts by fractional linear  transformations, and D7-branes serving as the source for the multi-valuedness of $\tau$. In this way, F-theory models correspond geometrically to torus fibrations over some compact base manifold. 
\par A well-known duality in string theory asserts that compactifying M-theory on a torus $\mathbf{T}^2$ with complex structure parameter $\tau$ and area $A$ is dual to the type IIB string compactified on a circle of radius $A^{-3/4}$ with axio-dilaton field $\tau$ \cites{MR1359982,MR1334520, MR1411455}.  This gives a connection between F-theory models and geometric compactifications of M-theory:  after compactifying an F-theory model further on $S^1$ without breaking supersymmetry, one obtains a model that is dual to M-theory compactified on the total space of the torus fibration. The geometric M-theory model preserves supersymmetry exactly when the total space of the family is a Calabi-Yau manifold.  In this way, we recover the familiar condition for supersymmetric F-theory models in eight dimensions or $D=8$: the total space of the fibration has to be a K3 surface.  
\par In this article, we will focus on F-theory models associated with eight-dimensional compactifications that correspond to genus-one fibrations with a section, or Jacobian elliptic fibrations, on algebraic K3 surfaces. As pointed out by Witten~\cite{MR1408164}, this subclass of models is in fact physically easier to treat since the existence of a section also implies the absence of NS-NS and R-R fluxes in F-theory.  Geometrically, the restriction to Jacobian elliptic fibrations facilitates model building with various non-Abelian gauge symmetries using the Tate algorithm \cites{MR1412112, MR2876044, MR3065908} where insertions of seven-branes in an F-theory model correspond to singular fibers in the M-theory model.  Through work of Kodaira  \cite{MR0184257} and N\'eron \cite{MR0179172}, all possible singular fibers in one-parameter families of elliptic curves have been classified. The catalog and its physical interpretation are by now well-known; see~\cite{MR3366121}.  As we shall see, our investigation of F-theory/heterotic string duality will be greatly aided by the existence of a fibration with section: it will allows us to utilize the mathematical classification of elliptic fibrations with section obtained by the authors in \cite{CHM19} and construct the dual non-geometric heterotic vacua in $D=8$ with two Wilson line parameters.
\par Our construction of F-theory models relies on the concrete relationship between modular forms on the moduli space of certain K3 surfaces of Picard rank $16$ and the equations of those K3 surfaces which have also been studied in \cites{MR1204828, MR4015343, Hosono_2019, CHM19}.  The K3 surfaces in question have a large collection of algebraic curve classes on them, generating a lattice known as $H \oplus E_7(-1) \oplus E_7(-1)$.  The presence of these classes restricts the form of moduli space, which turns out to be a space admitting {\em modular forms.}  The K3 surfaces turn out to be closely related to the family of double sextic surfaces, i.e., K3 surfaces obtained as double cover of the projective plane branched on a reducible sextic. The modular forms are therefore generalizations of the Siegel modular forms of genus two and can be constructed explicitly using the exceptional analytic equivalence between the bounded symmetric domains of type $IV_4$ and of type $I_{2,2}$. In fact, generators for the ring of modular forms were constructed by two of the authors in \cite{MR4015343}. In this article we construct explicit Weierstrass model for all inequivalent Jacobian elliptic fibrations supported on the family of K3 surfaces with $H \oplus E_7(-1) \oplus E_7(-1)$ lattice polarization and express their  parameters in terms of such modular forms.
\par This article is structured as follows: in Section~\ref{sec:family_of_K3s}  we introduce the fundamental mathematical object that is later used to describe (the dual of) the non-geometric heterotic string vacua in $D=8$ with two Wilson lines, namely the family of K3 surfaces with canonical $H \oplus E_7 (-1) \oplus E_7 (-1)$ polarization. We will show that this family admits exactly four inequivalent Jacobian elliptic fibrations. For each fibration we will construct a Weierstrass model whose parameters are modular forms generalizing well-known Siegel modular forms. In Section~\ref{sec:confluences} we present possible confluences of the singular fibers that can occur in the four Jacobian elliptic fibrations. Establishing these confluences turns out to be critical for matching the backgrounds in $D=8$ with the backgrounds dual to the heterotic string with an unbroken gauge algebra $\mathfrak{e}_8 \oplus \mathfrak{e}_8$ or $\mathfrak{so}(32)$ as we turn off the Wilson line parameters. In Section~\ref{sec:heterotic_models} we classify all non-geometric heterotic models obtained by the partial Higgsing (using two Wilson lines) of the heterotic gauge algebra $\mathfrak{g}=\mathfrak{e}_8 \oplus \mathfrak{e}_8$ or $\mathfrak{g}=\mathfrak{so}(32)$ for the associated low energy effective eight-dimensional supergravity theory, dual to the K3 surfaces from Section~\ref{sec:family_of_K3s}. There we find a surprise: as opposed to the Higgsing of the heterotic gauge algebra using only one Wilson line, the Higgsing with two Wilson lines produces two different branches for each type of heterotic string. We conclude the paper with a discussion of this surprise and its implications.
\section{A special family of K3 surfaces}
\label{sec:family_of_K3s}
Let $\mathcal{X}$ be a smooth complex algebraic K3 surface.   The group of divisors (modulo algebraic equivalence) is called the Ner\'on-Severi lattice of $\mathcal{X}$, denoted by $\operatorname{NS}(\mathcal{X})$. It is well known that $\operatorname{NS}(\mathcal{X})$ is an even lattice of signature $(1, p_\mathcal{X})$, where $p_\mathcal{X}$ is the Picard rank of $\mathcal{X}$ with $1 \le p_\mathcal{X} \le 20$. Following \cites{MR0357410,MR0429917,MR544937,MR525944,MR728992}, for a fixed even lattice $\mathrm{N}$ fo signature $(1, r)$, with $0 \le r\le 19$, we say that $\mathcal{X}$ is polarized by the lattice $\mathrm{N}$, if $\imath \colon \mathrm{N} \to \operatorname{NS}(\mathcal{X})$ is a primitive embedding of lattice for which $\imath(\mathrm{N})$ contains a pseudo-ample divisor class.   We call $(\mathcal{X}, \imath)$ an $\mathrm{N}$-polarized K3 surface.  Two $\mathrm{N}$-polarized K3 surfaces $(\mathcal{X}, \imath)$ and $(\mathcal{X}', \imath')$ are said to be isomorphic, if there exists an analytic isomorphism $\alpha: \mathcal{X} \to \mathcal{X}'$ such that $\alpha^* \circ \imath' = \imath$ where $\alpha^*$ is the appropriate morphism at cohomology level.  
\par This article aims to describe certain backgrounds for the heterotic string using a special class of such objects, namely,  K3 surfaces which are polarized by the rank sixteen lattice
\beq
\mathrm{N} = H \oplus E_7 (-1) \oplus E_7 (-1),
\eeq
where $H$ is the standard hyperbolic lattice of rank two (hyperbolic plane), and $E_7(-1)$ is the negative definite even lattice associated with the $E_7$ root system. K3 surfaces of this type are explicitly constructible: let $(\alpha, \beta, \gamma, \delta , \varepsilon, \zeta) \in \mathbb{C}^6 $ be a set of parameters, consider the projective quartic surface $\mathrm{Q}(\alpha, \beta, \gamma, \delta, \varepsilon, \zeta)$ in $\mathbb{P}^3(\mathbf{X}, \mathbf{Y}, \mathbf{Z}, \mathbf{W})$ defined by the homogeneous equation:   
\beq
\label{quartic1}
	\mathbf{Y}^2 \mathbf{Z} \mathbf{W}-4 \mathbf{X}^3 \mathbf{Z}+3 \alpha \mathbf{X} \mathbf{Z} \mathbf{W}^2+\beta \mathbf{Z} \mathbf{W}^3+\gamma \mathbf{X} \mathbf{Z}^2 \mathbf{W}- \frac{1}{2} \left (\delta \mathbf{Z}^2 \mathbf{W}^2+ \zeta \mathbf{W}^4 \right )+ \varepsilon \mathbf{X} \mathbf{W}^3 = 0 \,.
\eeq
The family in Equation~(\ref{quartic1}) was first introduced by Clingher and Doran in \cite{MR2824841} as a generalization of the Inose quartic in \cite{MR578868}.  Assuming that $(\gamma, \delta) \neq (0,0)$ and $ (\varepsilon, \zeta) \neq (0,0)$, it was proved in \cite{MR4015343} that the surface ${\mathcal{X}}(\alpha, \beta, \gamma, \delta , \varepsilon, \zeta)$ obtained as the minimal resolution of $(\ref{quartic1})$ is a K3 surface endowed with a canonical $\mathrm{N}$ polarization. All $\mathrm{N}$-polarized K3 surfaces, up to isomorphism, are in fact realized in this way. Moreover, one can tell precisely when two members of the above family are isomorphic.  Let $\mathfrak{G}$ be the subgroup of $\operatorname{Aut}(\mathbb{C}^6)$ generated by the set of transformations given below:
\beq
\label{eqn:action}
 (\alpha,\beta, \gamma, \delta, \varepsilon, \zeta) \ \longrightarrow \
(t^2 \alpha,  \ t^3 \beta,  \ t^5 \gamma,  \ t^6 \delta,  \ t^{-1} \varepsilon,  \ \zeta  ), \ {\rm with} \ t \in \mathbb{C}^* $$
$$ (\alpha,\beta, \gamma, \delta, \varepsilon, \zeta) \ \longrightarrow \  
(\alpha,  \beta,  \varepsilon,  \zeta, \gamma,  \delta  ) \,.
\eeq
It then follows that two K3 surfaces in the above family are isomorphic if and only if their six-parameter coefficient sets belong to the same orbit of $\mathbb{C}^6$ under $\mathfrak{G}$; see  \cite{MR4015343} .
\subsection{Jacobian elliptic fibrations on \texorpdfstring{$\mathcal{X}$}{X}} 
\label{ssec:JAC-fibrations}
Recall that a Jacobian elliptic fibration on $\mathcal{X}$ is a pair $(\pi, \sigma)$ consisting of a proper map of analytic spaces $\pi : \mathcal{X} \to \mathbb{P}^1$, whose generic fiber is a smooth genus one curve, and a section $\sigma : \mathbb{P}^1 \to \mathcal{X}$ in the elliptic fibration $\pi$. If $\sigma'$ is another section of the Jacobian fibration $(\pi, \sigma)$ then there exists an automorphism of $\mathcal{X}$ that preserves $\pi$ and maps $\sigma$ to $\sigma'$.  By identifying the set of sections of $\pi$ and the group of automorphisms of $\mathcal{X}$ preserving $\pi$, the set of all sections form a group, known as the Mordell-Weil group of the Jacobian elliptic fibration, denoted by $\mathrm{MW}(\pi, \sigma)$.  
\par Classifying Jacobian elliptic fibrations on $\mathcal{X}$ corresponds to classifying primitive lattice embeddings $H \hookrightarrow \operatorname{NS}(\mathcal{X})$ since isomorphism classes of Jacobian elliptic fibrations on $\mathcal{X}$ are in one-to-one correspondence with isomorphism classes of primitive lattice embeddings $H \hookrightarrow \operatorname{NS}(\mathcal{X})$ \cite{MR2369941}*{Lemma 3.8}.  A lattice theoretic analysis by the authors revealed that there are exactly four such (non-isomorphic) primitive lattice embeddings \cite{CHM19}*{Prop.~2.3}.  Consequently, it was shown that an $\mathrm{N}$-polarized K3 surface carries four Jacobian elliptic fibrations, up to automorphisms \cite{CHM19}*{Thm.~3.6}. In Sections~\ref{sec:StdFibration}-\ref{sec:MaxFibration} we will briefly review the construction of Weierstrass models for these fibrations obtained in \cite{CHM19}*{Thm.~3.6}. We use the Kodaira classification of singular fibers to describe the four Jacobian elliptic fibrations \cite{MR0184257}, which we call the standard, alternate, base-fiber dual, and maximal fibration.
\subsubsection{The standard fibration.} \label{sec:StdFibration}
Substituting
\beq
\mathbf{X} = u v x\,, \quad  \mathbf{Y}= y\,, \quad \mathbf{Z} = 4 u^4 v^2 z\,, \quad \mathbf{W} = 4 u^3 v^3 z \,,
\eeq
in Equation~(\ref{quartic1}) yields the Jacobian elliptic fibration $\pi_{\text{std}} : \mathcal{X} \rightarrow \mathbb{P}^1$ with fiber $\mathcal{X}_{[u:v]}$, given by the Weierstrass equation 
\beq
\label{eqn:std}
\mathcal{X}_{[u:v]}: \quad y^2 z  = x^3 +  f(u, v) \, x z^2 + g(u, v) \, z^3 \,, 
\eeq
equipped with the section $\sigma_{\text{std}} : [x:y:z] = [0:1:0]$, and with a discriminant 
\beq
\Delta(u, v) =  4 f^3 + 27 g^2 =  64 \, u^9 v^9  p(u, v) \,,
\eeq
where
\beq
f(u, v)  = - 4 u^3 v^3  \Big(\gamma  u^2 + 3  \alpha  uv + \varepsilon v^2\Big) \,, \quad
g(u, v)  =  8  u^5 v^5 \Big(\delta  u^2 - 2  \beta  u v  + \zeta  v^2\Big) \,,
\eeq
and $p(u, v) = 4  \gamma^3 u^6 + \dots + 4 \varepsilon^3 v^6$ is an irreducible homogeneous polynomial of degree six. 
\par Equation~(\ref{eqn:std}) defines a Jacobian elliptic fibration with six singular fibers of Kodaira type $I_1$, two singular fibers of Kodaira type $III^*$ (ADE type $E_7$), and a trivial Mordell-Weil group of sections $\operatorname{MW}(\pi_{\text{std}}, \sigma_{\text{std}}) = \{ \mathbb{I}\}$.
\subsubsection{The alternate fibration.}\label{sec:AltFibration}
Substituting
\beq
\mathbf{X}= 2 u v x\,, \quad  \mathbf{Y}=  y \,, \quad \mathbf{Z} = 4 v^5 (-2 \varepsilon u+\zeta v) z\,, \quad \mathbf{W} = 2 v^2 x\,,
\eeq
into Equation~(\ref{quartic1}), determines the Jacobian elliptic fibration  $\pi_{\text{alt}}:  \mathcal{X} \rightarrow \mathbb{P}^1$ with fiber $\mathcal{X}_{[u:v]}$, given by the equation 
\beq
\label{eqn:alt}
\mathcal{X}_{[u:v]}: \quad y^2 z  = x \Big(x^2 + A(u, v) \, x z +B(u, v) \, z^2 \Big) \,, 
\eeq
equipped with the section $\sigma_{\text{alt}} : [x:y:z] = [0:1:0]$, the two-torsion section  $[x:y:z] = [0:0:1]$, and with a discriminant 
\beq
\Delta(u, v) =  B(u, v)^2 \, \Big(A(u, v)^2-4 B(u, v)\Big) \,,
\eeq
where
\beq
\begin{split}
	A(u, v) = 4 v (4 u^3-3\alpha  uv^2- \beta v^3) \,,\quad
	B(u, v) = 4 v^6 (2\gamma u-\delta v)(2\varepsilon u -\zeta v)   \,.
\end{split}
\eeq
\par Equation~(\ref{eqn:alt}) defines a Jacobian elliptic fibration with six singular fibers of Kodaira type $I_1$,  two singular fibers of Kodaira type $I_2$ (ADE type $A_1$), and a singular fiber of Kodaira type $I_8^*$ (ADE type $D_{12}$), and a Mordell-Weil group of sections $\operatorname{MW}(\pi_{\text{alt}}, \sigma_{\text{alt}}) = \mathbb{Z}/2\mathbb{Z}$.
\subsubsection{The base-fiber-dual fibration}
\label{sec:BFDFibration}
Substituting
\beq
\begin{array}{lclclcl}
	\mathbf{X} &=& 3 uv (x+6\gamma \varepsilon uv^3z)\,, 							&\quad&	\mathbf{Y}&=&  y \,, \\ [0.2em]
	\mathbf{Z} &=& 6 v^2 (\varepsilon x-6\gamma \varepsilon^2 uv^3z-18 \zeta u^2 v^2 z)\,,	&\quad&	\mathbf{W}&=&108 u^3 v^3 z \,,
\end{array}
\eeq
into Equation~(\ref{quartic1}) determines a Jacobian elliptic fibration  $\pi_{\mathrm{bfd}}: \mathcal{X} \to \mathbb{P}^1$ with fiber $\mathcal{X}_{[u:v]}$, given by
the equation
\beq
\label{eqn:bfd}
\mathcal{X}_{[u:v]}: \quad y^2 z = x^3 + F(u, v) \, x z^2 + G(u, v) \, z^3 \,,
\eeq
admitting the section $\sigma_{\text{bfd}} : [x:y:z] = [0:1:0]$, and with a discriminant 
\beq
\Delta(u, v) =  4  F^3 + 27  G^2 =  - 2^6 3^{12} u^8 v^{10}  P(u, v) \,,
\eeq
where
\beq
\begin{split}
	F(u, v)  = & \; - 108\, u^2 v^4  \Big(9 \alpha  u^2 - 3 (\gamma \zeta+\delta \varepsilon)uv + \gamma^2 \varepsilon^2 v^2\Big) \,, \\
	G(u, v) = & \; - 216 u^3 v^5 \Big( 27 u^4 + 54 \beta u^3 v + 27(\alpha \gamma \varepsilon+ \delta \zeta) u^2 v^2 \\
	& \; - 9 \gamma \varepsilon (\gamma\zeta+\delta\varepsilon) u v^3 + 2 \gamma^3 \varepsilon^3 v^4 \Big) \,,
\end{split}
\eeq
and $P(u, v) =   \gamma^2\varepsilon^2 (\gamma\zeta-\delta\varepsilon)^2 v^6 + O(u)$ is an irreducible homogeneous polynomial of degree six.
\par Equation~(\ref{eqn:bfd}) defines a Jacobian elliptic fibration with six singular fibers of Kodaira type $I_1$,  one singular fibre of Kodaira type $I_2^*$ (ADE type $D_6$), and a singular fiber of Kodaira type $II^*$ (ADE type $E_8$), and a Mordell-Weil group of sections $\operatorname{MW}(\pi_{\text{bfd}}, \sigma_{\text{bfd}}) =  \{\mathbb{I}\}$.
\subsubsection{The maximal fibration} \label{sec:MaxFibration}
The maximal fibration is induced by intersecting the quartic $\mathrm{Q}(\alpha, \beta, \gamma, \delta, \varepsilon, \zeta)$ with the pencil of quadric surfaces 
\beq 
\label{eq:conic}
\begin{split}
	\mathrm{C}_3(u, v)= v \Big( 2 \gamma^2 \delta \varepsilon \zeta \mathbf{X} \mathbf{Z} + (6\alpha\gamma\delta\varepsilon\zeta +4 \beta \gamma \delta \varepsilon^2  + 4 \beta \gamma^2 \varepsilon \zeta+2 \delta^2 \zeta^2) \mathbf{X} \mathbf{W} - \gamma\delta^2\varepsilon \zeta \mathbf{Z} \mathbf{W} \\
	+ 2 \gamma \delta \varepsilon \zeta \mathbf{Y}^2  - (8 \beta \gamma^2 \varepsilon^2 + 4 \delta^2 \varepsilon \zeta + 4 \gamma \delta \zeta^2) \mathbf{X}^2 \Big) 
	+ u (2 \gamma \mathbf{X} - \delta \mathbf{W})(2 \varepsilon \mathbf{X} - \zeta \mathbf{W}) =0 \,,
\end{split}
\eeq
with $[u : v] \in \mathbb{P}^1$.   Making the substitutions 
\beq
\begin{split}
	\mathbf{X}= \delta\zeta v \big((2\beta \gamma \varepsilon v-u)x -2 \gamma \delta^5 \varepsilon \zeta^5 v^5 z\big)\,, \quad  \mathbf{Y}=  y \,, \quad \mathbf{W} =  2 \delta^2\zeta^2 v^2 x \,,
\end{split}
\eeq
and $\mathbf{Z}=\mathbf{Z}(x, y, z, u, v)$, obtained by solving Equation~(\ref{eq:conic}) for $\mathbf{Z}$, determines a Jacobian elliptic fibration  $\pi_{\mathrm{max}}: \mathcal{X} \to \mathbb{P}^1$ with fiber $\mathcal{X}_{[u:v]}$, given by the equation
\beq
\label{eqn:max}
\mathcal{X}_{[u:v]}: \quad y^2 z = x^3 + a(u, v) \, x^2 z + b(u, v) \, x z^2 + c(u, v) \, z^3 \,,
\eeq
admitting the section $\sigma_{\text{max}} : [x:y:z] = [0:1:0]$, and with the discriminant 
\beq
\begin{split}
	\Delta(u, v) & =  b^2 \big(a^2 - 4 b\big) - 2 a c  \big(2 a^2-9 b \big) -27 c^2  = 64 \delta^{16} \zeta^{16} v^{16} d(u, v)\,,
\end{split}
\eeq
where
\beq
\begin{split}
	a(u, v)  & = - 2 \delta \zeta v \Big( u^3 - 6 \beta \gamma \varepsilon u^2 v + 3(4 \beta^2\gamma^2 \varepsilon^2- \alpha\delta^2\zeta^2) uv^2\\
	& \quad - 2 \beta ( 4 \beta^2\gamma^3 \varepsilon^3 - 3 \alpha\gamma \delta^2 \varepsilon\zeta^2-\delta^3\zeta^3) v^3 \Big)\,, \\
	b(u, v)  & = - 4 \delta^6 \zeta^6 v^6 \Big( 2 \gamma \varepsilon u^2 - (8 \beta \gamma^2 \varepsilon^2+ \gamma\delta\zeta^2+\delta^2\varepsilon\zeta) uv \\
	& \quad +(8 \beta^2 \gamma^3 \varepsilon^3-3 \alpha\gamma\delta^2\varepsilon\zeta^2 +2  \beta \gamma^2 \delta \varepsilon \zeta^2 + 2 \beta\gamma \delta^2\varepsilon^2\zeta -\delta^3\zeta^3)v^2 \Big) \,, \\
	c(u, v) &= - 8 \gamma \delta^{11} \varepsilon \zeta^{11} v^{11} \Big( \gamma \varepsilon u - (2 \beta \gamma^2 \varepsilon^2 + \gamma\delta \zeta^2+\delta^2\varepsilon\zeta)v\Big) \,,
\end{split}
\eeq
and $d(u, v) =  (\gamma\zeta - \delta\varepsilon)^2 u^8 + O(v)$ is an irreducible homogeneous polynomial of degree eight.
\par Equation~(\ref{eqn:max}) defines a Jacobian elliptic fibration with eight singular fibers of Kodaira type $I_1$,  one singular fiber of Kodaira type $I_{10}^*$ (ADE type $D_{14}$), and a Mordell-Weil group of sections $\operatorname{MW}(\pi_{\text{max}}, \sigma_{\text{max}}) =  \{\mathbb{I}\}$.
\subsection{Modular Description}
\label{ssec:modular}
The parameters of the defining equations for the Weierstrass models in Sections~\ref{sec:StdFibration}-\ref{sec:MaxFibration} can be interpreted as modular forms, established by two of the authors and Shaska in~\cite{MR4015343}. 
\par Let $L^{2,4}$ be the orthogonal complement $\mathrm{N}^\perp \subset \Lambda_{K3}$ in the K3 lattice $\Lambda_{K3}= H^{\oplus 3} \oplus E_8(-1) \oplus E_8(-1)$ with orthogonal transformations $O(L^{2,4})$. Let $\mathcal{D}_{2,4}$ be the Hermitian symmetric space, specifically the bounded symmetric domain of type $IV_4$, given as
\beq
 \mathcal{D}_{2,4} = O^+(2,4) / \big(  SO(2) \times  O(4) \big) \,,
\eeq  
where $O^+(2,4)$ denotes the subgroup of index two of the pseudo-orthogonal group $O(2,4)$ consisting of the elements whose upper left minor of order two is positive. Let $O^+(2,4;\mathbb{Z})=O(L^{2,4}) \cap O^+(2,4)$ be the arithmetic lattice of $O^+(2,4)$, i.e., the discrete cofinite group of holomorphic automorphisms on the bounded Hermitian symmetric domain $\mathcal{D}_{2,4}$. We also set $SO^+(2,4) =O^+(2,4) \cap SO(2,4)$ and $SO^+(2,4;\mathbb{Z})= O(L^{2,4}) \cap SO^+(2,4)$.

 An appropriate version of the Torelli theorem \cite{MR1420220} gives rise to an analytic isomorphism between the moduli space of $\mathrm{N}$-polarized K3 surfaces and the quasi-projective four-dimensional algebraic variety $\mathcal{D}_{2,4}/O^+(2,4;\mathbb{Z})$. We consider the normal, finitely generated algebra $A(\mathcal{D}_{2,4},\mathcal{G}) = \oplus_{k \ge 0} A(\mathcal{D}_{2,4},\mathcal{G})_k$ of automorphic forms on $\mathcal{D}_{2,4}$ relative to a discrete subgroup $\mathcal{G}$ of finite covolume in $O^+(2,4)$, graded by the weight $k$ of the automorphic forms.  For $\mathcal{G}= O^+(2,4;\mathbb{Z})$, the algebra $A(\mathcal{D}_{2,4},\mathcal{G})$ is freely generated by forms $J_k$ of weight $2k$ with $k=2, 3, 4, 5, 6$; this is a special case of a general result proven by Vinberg in \cites{MR2682724, MR3235787}. A subtle point here is that one has to obtain $A(\mathcal{D}_{2,4}, \mathcal{G})$ as the \emph{even part}
\beq
\label{eq:even_part}
 A(\mathcal{D}_{2,4}, \mathcal{G}) = \Big[  A(\mathcal{D}_{2,4}, \mathcal{G}_0) \Big]_{\text{even}}
\eeq
of the ring of automorphic forms with respect to subgroup $\mathcal{G}_0 = SO^+(2,4;\mathbb{Z})$ of index two.

\par Based on the exceptional analytic equivalence between the bounded symmetric domains of type $IV_4$ and of type $I_{2,2}$, an explicit description of the generators $\{J_k\}_{k=2}^6$ can be derived. To start, we remark that $\mathcal{D}_{2,4} \cong \mathbf{H}_{2,2}$, where 
\beq
\label{Siegel_tau}
\mathbf{H}_{2,2} =  \left \lbrace \begin{pmatrix}
	\tau_1 & z_1\\z_2 & \tau_2
	\end{pmatrix} 
	\in \text{Mat}(2, 2; \mathbb{C}) \; \Big| \; \operatorname{Im}{(\tau_1)} \cdot \operatorname{Im}{(\tau_2)} > \frac{1}{4} |z_1 - \bar{z}_2|^2, \; \operatorname{Im}{\tau_1}> 0 \right\rbrace \,.
\eeq
The domain $\mathbf{H}_{2,2}$ is a generalization of the Siegel upper-half space $\mathbb{H}_2$ in the sense that 
\beq
\label{eqn:restriction}
\mathbb{H}_2 = \left \lbrace \varpi \in \mathbf{H}_{2,2} \; \big| \; \varpi^t = \varpi \right\rbrace \,.
\eeq
A subgroup $\Gamma \subset \operatorname{U}(2,2)$, given by
\beq
\label{modular_group}
\Gamma =  \left\lbrace G \in \operatorname{GL}\big(4,\mathbb{Z}[i]\big) \,  \Big| \, G^\dagger \cdot \left( \begin{array}{cc} 0 & \mathbb{I}_2 \\  -\mathbb{I}_2 & 0 \end{array}\right) \cdot G = \left( \begin{array}{cc} 0 & \mathbb{I}_2 \\  -\mathbb{I}_2 & 0 \end{array}\right)  \right\rbrace \,,
\eeq
acts on $\varpi \in \mathbf{H}_{2,2}$ by 
\beq
\forall \, G= \left(\begin{matrix} A & B \\ C & D \end{matrix} \right) \in \Gamma: \quad G\cdot \varpi=(C\cdot \varpi+D)^{-1} (A\cdot \varpi+B)  \,.
\eeq
There is an involution $T$ acting on $\mathbf{H}_{2,2}$ by transposition, i.e., $\varpi \mapsto T\cdot \varpi=\varpi^t$, yielding an extended group as the semi-direct product $\Gamma_T = \Gamma \rtimes \langle T \rangle$. Moreover, the group $\Gamma_T$ has the index-two subgroup given by
\beq
\label{index-two}
\begin{split}
 \Gamma^{+}_T =   \Big\lbrace g = G\, T^{n} \in  \Gamma_T  & \ \Big| \ n \in \lbrace 0,1 \rbrace, \; (-1)^n \, \det{G} =1 \Big\rbrace \,.
 \end{split}
\eeq
The identification $\mathcal{D}_{2,4} \cong \mathbf{H}_{2,2}$ gives rise to  an homomorphism $\operatorname{U}(2,2)  \to SO^+(2,4)$ which also identifies $\mathcal{G}_0 \cong  \Gamma^{+}_T$ \cite{MR2682724}.
\par In the above context, the five modular forms $\{J_k\}_{k=2}^6$ can be computed in terms of theta functions on $\mathbf{H}_{2,2}$ introduced by Matsumoto~\cite{MR1204828}. The result is that the generators of the graded ring of even modular forms relative to the group $\Gamma_T^+$ can be determined explicitly:
\begin{theorem}[\cite{MR4015343}]
\label{thm:ROMF}
The invariants  $J_2,  J_3, J_4, J_5, J_6$ are modular forms of respective weights 4, 6, 8, 10, and 12, relative to the group $\Gamma_T^+$. 
\end{theorem}
\noindent
We remark that there also is an automorphic form $\mathfrak{a}$ with nontrivial automorphic factor of weight $10$ satisfying $\mathfrak{a}^2 =J_5^2 -4  J_4 J_6$; for its precise definition we refer to \cite{MR4015343}.
\par The fact that two K3 surfaces in the family in Equation~(\ref{quartic1}) are isomorphic if and only if their six-parameter coefficient sets belong to the same orbit of $\mathbb{C}^6$ under $\mathfrak{G}$ in Equation~(\ref{eqn:action}) allows to identify the modular forms $J_{k}$ relative to the action of $\Gamma_T$ with the following set of invariants associated to the K3 surfaces in the family:
\beq
\label{modinv}
J_2 = \alpha, \ \ \ J_3 = \beta, \ \ \ J_4 =  \gamma \cdot \varepsilon, \ \ \ J_5 = \gamma \cdot \zeta + \delta \cdot \varepsilon, \ \ \  J_6 = \delta \cdot \zeta \,,
\eeq
which then allows one to prove:
\begin{theorem}[\cite{MR4015343}]
 The four-dimensional open analytic space
\beq
\label{modulispace}
\mathfrak{M}_\mathrm{N} \ = \ \Big \{ \ 
\vec{J} = \left [ \ J_2 , \ J_3 , \  J_4, \ J_5 , \ J_6  \ \right ]   \in  \mathbb{W}\mathbb{P}(2,3,4,5,6) \ \vert  \ 
( J_4 , \;  J_5, \; J_6 ) \neq (0,0,0) \ 
\Big \} \ 
\eeq
forms a coarse moduli space for $\mathrm{N}$-polarized K3 surfaces. 
\end{theorem}
\par Under the restriction of $\mathbf{H}_{2,2} / \Gamma^+_{T}$ to $\mathbb{H}_2 / \text{Sp}_4(\mathbb{Z})$ induced by Equation~(\ref{eqn:restriction}) we obtain
\beq
   \label{SiegelRestriction}
  [J_2(\varpi) : J_3(\varpi) : J_4(\varpi) : J_5(\varpi) : J_6(\varpi)] =
  [\psi_4(\tau) : \psi_6(\tau) : 0 : 2^{12} 3^5 \chi_{10}(\tau) : 2^{12} 3^6 \chi_{12}(\tau)], 
\eeq
where $\psi_4, \psi_6, \chi_{10},$ and $\chi_{12}$ are the Siegel modular forms of respective weights $4, 6, 10$ and $12$ introduced and defined by Igusa in \cite{MR0141643}.
\par For the four inequivalent Jacobian elliptic fibrations on the family of N-polarized K3 surfaces $\mathcal{X}$ obtained in Section~\ref{ssec:JAC-fibrations}, we will now construct Weierstrass models with coefficients in $\mathbb{Q}[J_2, J_3, J_4, J_5, J_6]$ or $\mathbb{Q}[J_2, J_3, \mathfrak{a}, J_5, J_6]$ in the case of the standard fibration.  
\subsubsection{The standard fibration}\label{sec:StdFibMF} 
Assuming $J_6 \ne 0$ we denote the two solutions of the equation $\mathfrak{a}^2 =J_5^2 -4  J_4 J_6$ by $\pm \mathfrak{a}$. The elliptic fibration  $\pi_{\mathrm{std}} \colon \mathcal{X} \to \mathbb{P}^1$ in Section~\ref{sec:StdFibration} can be written in a suitable affine coordinate chart as
\beq
\label{eqn:std_mod}
Y^2 = X^3  + f_\pm (t) \, X + g(t)  \,, 
\eeq
with
\beq
f_\pm(t)  = -  t^3 J_6^3 \left(\frac{J_5 \mp \mathfrak{a}}{2}  t^2 + 3 J_2 J_6   t + \frac{(J_5 \pm \mathfrak{a}) J_6}{2} \right) \, , \quad 
g(t) =   J_6^5 t^5  \Big(t^2 - 2 J_3 t + J_6 \Big) \, ,
\eeq
and a discriminant $\Delta = J_6^{9} t^9 p_\pm(t)$  where 
\beq
p_\pm(t) = 2 \Big( J_5^2 (J_5 \pm \mathfrak{a}) - J_4 J_6 (3J_5 \pm \mathfrak{a})   \Big) t^6 + \dots + 2 J_6^3\Big( J_5^2 (J_5 \mp \mathfrak{a}) -J_4 J_6 (3J_5 \mp  \mathfrak{a}) \Big) \,.
\eeq
Notice that we have $f_-(J_6/t)= J_6^4 f_+(t)/t^8$ and $g(J_6/t) = J_6^6 g(t)/t^{12}$. Since the map 
\beq
  (t, X, Y) \mapsto (t', X', Y') = (J_6/t \ , \, J_6^2 X/t^4 \ , \ -J_6^3 Y/t^6) \,,
\eeq
maps the K3 surfaces in Equation~(\ref{eqn:std_mod}) with $f_-(t)$ to one with $f_+(t)$ and the holomorphic two form $dt \wedge dX/Y$ to $dt' \wedge dX'/Y'$, it provides a holomorphic, symplectic morphism between the two K3 surfaces. 
\par For $J_6=0$, we have $\delta=0$ or $\zeta=0$, and the elliptic fibration  $\pi_{\mathrm{std}} \colon \mathcal{X} \to \mathbb{P}^1$ in Section~\ref{sec:StdFibration}  can be written as either
\beq
Y^2 = X^3  - t^3 \Big(t^2 + 3 J_2 \,t + J_4 \Big) \, X + t^5 \Big( J_5 - 2 J_3  t \Big) \,,
\eeq
or
\beq
Y^2 = X^3  - t^3 \Big(J_4 t^2 + 3 J_2 \,t + 1 \Big) \, X + t^5 \Big( - 2 J_3  t + J_5  t^2\Big) \,.
\eeq 
The fibrations are related by the birational morphism
\beq
  (t, X, Y) \mapsto (t', X', Y') = (1/t \ , \, X/t^4 \ , \ -Y/t^6) \,,
\eeq
which also maps the holomorphic two-form $dt \wedge dX/Y$ to $dt' \wedge dX'/Y'$. 
\subsubsection{The alternate fibration}\label{sec:AltFibMF} 
The Jacobian elliptic fibration  $\pi_{\mathrm{alt}}: \mathcal{X} \to \mathbb{P}^1$ in Section~\ref{sec:AltFibration} is written in a suitable affine coordinate chart as
\beq
\label{eqn:alt_mod}
Y^2  = X \Big( X^2  + A(t) \,  X + B(t) \Big)  \,,
\eeq
with
\beq
A(t) = t^3-3 J_2 \, t-2  J_3 \,, \quad \quad B(t) = J_4 \, t^2-J_5 \, t + J_6 \,,
\eeq 
and a discriminant $\Delta = E(t)^2 D(t)$ where $E(t)=J_4 t^2 -J_5 t + J_6$ and
\beq
D(t) = t^6 - 6 J_2 t^4 - 4 J_3 t^3 + (9J_2^2-4J_4) t^2 + (12 J_2 J_3+4 J_5) t + 4 (J_3^2-J_6) \,.
\eeq
\subsubsection{The base-fiber-dual fibration}\label{sec:BFDFibMF} 
The Jacobian elliptic fibration  $\pi_{\mathrm{bfd}}: \mathcal{X} \to \mathbb{P}^1$ in Section~\ref{sec:BFDFibration} is written in a suitable affine coordinate chart as
\beq
\label{eqn:bfd_mod}
Y^2 = X^3 + F(t) \, X + G(t) \,,
\eeq
with
\beq
\begin{split}
	F(t)  = & \; t^2 \Big( - 3 J_2 t^2 - J_5 t -\frac{1}{3} J_4^2 \Big)\,, \\
	G(t) = & \; t^3 \Big( t^4 -2 J_3 t^3 + (J_2J_4+J_6) t^2 + \frac{1}{3} J_4J_5 t+ \frac{2}{27} J_4^3 \Big) \,,
\end{split}
\eeq
and a discriminant $\Delta = t^8 P(t)$ where $P(t) =  -27 t^6 +108  J_3 t^5 + \dots + \mathfrak{a}^2 J_4^2$.
\subsubsection{The maximal fibration}\label{sec:MaxFibMF} 
The Jacobian elliptic fibration  $\pi_{\mathrm{max}}: \mathcal{X} \to \mathbb{P}^1$ in Section~\ref{sec:MaxFibration} is written in a suitable affine coordinate chart as
\beq
\label{eqn:max_mod}
Y^2  = X^3 + a(t) \, X^2 + b(t) \, X + c(t)  \,,
\eeq
with
\beqn
\begin{split}
	a(t)  & = J_6 \Big( t^3 + 6 J_3 J_4 t^2 + 3  ( 4 J_3^2J_4^2 - J_2 J_6^2) t - 2 J_3 (3 J_2 J_4 J_6^2-4J_3^2J_4^3+J_6^3) \Big)\,, \\
	b(t)  & = - J_6^6 \Big( 2 J_4  t^2 + (8 J_3J_4^2+J_5 J_6) t + (8J_3^2J_4^3-3J_2J_4J_6^2+2 J_3J_4J_5J_6-J_6^3)\Big) \,, \\
	c(t) &=  J_4 J_6^{11} \Big( J_4 t + (2J_3J_4^2+J_5J_6)\Big) \,,
\end{split}
\eeqn
and a discriminant $\Delta = J_6^{16} d(t)$ where $d(t) =  \mathfrak{a}^2 t^8 + \dots$ is an irreducible polynomial of degree eight. 
By an appropriate change of coordinates, one can write the fibration in Equation~(\ref{eqn:max_mod}) in Weierstrass normal form
\beqn
y^2 = x^3 + \alpha(t) x + \beta(t) 
\eeqn
where $\alpha(t) = {t}^{6}+ \dots$ and $\beta(t) = {t}^{9}-{\dots}
$ are irreducible polynomials of degree six and nine respectively.

A simple computation shows that the various discriminants (denoted by $\operatorname{Disc}_t$) and resultants (denoted by $\operatorname{Res}_t$) with respect to the variable $t$ are related to a modular form $J_{30}$ of weight 60 which is a polynomial in $\{J_k\}_{k=2}^6$ and given by 
\beq
\label{eqn:J_30}
J_{30}= \operatorname{Disc}_t D = \operatorname{Disc}_t d = \frac{2^4}{3^{18} J_6^{30}} \frac{\operatorname{Disc}_t p}{\operatorname{Res}_t^3(t^{-3}f, \, t^{-5}g)} = -\frac{J_2^9}{3^{21}}\frac{\operatorname{Disc}_t P}{\operatorname{Res}_t^3(t^{-2}F, \, t^{-3}G)}   \,.
\eeq
\section{Confluences of Singular Fibers}
 \label{sec:confluences}
In this section we describe some confluences of the singular fibers that are possible for the fibrations constructed in Sections \ref{sec:StdFibration} -\ref{sec:MaxFibration}. The results are summarized in Table~\ref{fig:polarization}. In the table, $D(\Lambda)$ denotes the discriminant group of a lattice $\Lambda$ which is an important group theoretic invariant of the lattice used in Nikulin's classification theory \cites{MR544937,MR633160}.
\par For $J_4=0$ the \emph{base-fiber-dual fibration} in Equation~(\ref{eqn:bfd_mod}) specializes to 
\beq
\label{eqn:limit_std}
Y^2 = X^3  - t^3 \Big( 3 J_2 \,t + J_5 \Big) \, X + t^5 \Big(  t^2 - 2 J_3  t + J_6 \Big) \,.
\eeq 
For $J_4=0$, $J_6 \not =0$, we have $\mathfrak{a}= \pm J_5$: the \emph{standard fibration} in Equation~(\ref{eqn:std_mod}) then simplifies, after rescaling, to either Equation~(\ref{eqn:limit_std}) or
\beq
Y^2 = X^3  - t^4 \Big( J_5 t + 3 J_2 \,t \Big) \, X + t^5 \Big( J_6 t^2 - 2 J_3  t + 1 \Big) \,.
\eeq
The two fibrations are related by the birational morphism
\beq
  (t, X, Y) \mapsto (t', X', Y') = (1/t \ , \, X/t^4 \ , \ -Y/t^6) \,,
\eeq
which also maps the holomorphic two-form $dt \wedge dX/Y$ to $dt' \wedge dX'/Y'$.  Thus, the standard and the base-fiber-dual fibration specialize to the same elliptic fibration in Picard rank $17$ ($J_4=0$) and  Picard rank $18$ ($J_4=J_5=0$), marked in color-coded rows (blue) in Table~\ref{fig:polarization}. These are precisely the equations derived in \cite{MR3366121} for the F-theory dual of a non-geometric heterotic theory with gauge algebra $\mathfrak{g} = \mathfrak{e}_8 \oplus \mathfrak{e}_7$ and $\mathfrak{g} = \mathfrak{e}_8 \oplus \mathfrak{e}_8$ one non-vanishing Wilson line parameter. We will explain the physical relevance of this observation in Section~\ref{sec:heterotic_models}.
\par For $J_4=0$ the \emph{alternate fibration} in Equation~(\ref{eqn:alt_mod}) specializes to 
\beq
\label{eqn:limit_alt}
Y^2 = X^3  + \Big( t^3 - 3 J_2 \,t - 2 J_3 \Big) \, X^2 - \Big(  J_5 t - J_6 \Big) \,.
\eeq 
Similarly, for $J_4=0$, $J_6 \not =0$ the \emph{maximal fibration} in Equation~(\ref{eqn:max_mod}) simplifies, after rescaling, to Equation~(\ref{eqn:limit_alt}). Thus, the alternate and the maximal fibration specialize to the same fibration in Picard rank $17$ ($J_4=0$) and  Picard rank $18$ ($J_4=J_5=0$), marked in color-coded rows (red) in Table~\ref{fig:polarization}. This is precisely the equation derived in \cite{MR3366121} for the F-theory dual of a non-geometric heterotic theory with gauge algebra $\mathfrak{g} = \mathfrak{e}_8 \oplus \mathfrak{e}_7$ and one non-vanishing Wilson line parameter. We will explain the physical relevance of this observation in Section~\ref{sec:heterotic_models}.
\par Moreover, Table~\ref{fig:polarization} gives the confluences of singular fibers that occur along the vanishing loci of $\mathfrak{a}=0$ (with $\mathfrak{a}^2 = J_5^2 - 4 J_4 J_6$) and $J_{30}=0$ (where $J_{30}$ is defined in Equation~(\ref{eqn:J_30})). We also determined the possible confluences $2 I_1 \to II$ and $I_2 + I_1 \to III$ that can occur within the four elliptic fibrations (the polynomials $f, g$, $D, E$, $F, G$, and $\alpha, \beta$ are defined in Sections \ref{sec:StdFibration} -\ref{sec:MaxFibration}).
 
\begin{figure}
	\begin{subfigure}{\textwidth}
		\centering
		\scalebox{\MyScale}{
			\begin{tabular}{| c || c  | c | c | c | c|}
				\hline 
				&&&&& \\[-0.9em]
				Fibration (\ref{eqn:std_mod}) & $p_\mathcal{X}$ & Singular Fibers & $\operatorname{MW}(\pi_{\text{std}}, \sigma_{\text{std}})$ 
				& Lattice Polarization $\Lambda$ & $D(\Lambda)$ \\[0.1em]
				\hline\hline	
				&&&&& \\[-0.9em]
				generic & 16 & $2 III^* + 6 I_1$ & $\{\mathbb{I}\}$
				& $H \oplus E_7(-1) \oplus E_7(-1)$ & $\mathbb{Z}_2^2$\\[0.1em]
				\hline
				&&&&& \\[-0.9em]
				$\operatorname{Res}_t(t^{-3}f, \, t^{-5}g)=0$  & 16 & $2 III^* + II + 4 I_1$ & $\{\mathbb{I}\}$
				& $H \oplus E_7(-1) \oplus E_7(-1)$ & $\mathbb{Z}_2^2$\\[0.1em]
				\hline
				&&&&& \\[-0.9em]
				$J_{30} =0$ & 17 & $2 III^* + I_2 + 4  I_1$ & $\{\mathbb{I}\}$
				& $H \oplus E_7(-1) \oplus E_7(-1)  \oplus A_1(-1)$ & $\mathbb{Z}_2^3$\\[0.1em]
				\hline
				&&\cellcolor{blue!25}&\cellcolor{blue!25}&\cellcolor{blue!25}&\cellcolor{blue!25} \\[-0.9em]
				$J_4=0$ & 17 & \cellcolor{blue!25}  $II^* + III^* + 5 I_1$ & \cellcolor{blue!25}  $\{\mathbb{I}\}$ 
				& \cellcolor{blue!25}  $H \oplus E_8(-1) \oplus E_7(-1)$ & \cellcolor{blue!25}  $\mathbb{Z}_2 $\\[0.1em]
				\hline
				&&\cellcolor{blue!25}&\cellcolor{blue!25}&\cellcolor{blue!25}&\cellcolor{blue!25} \\[-0.9em]
				$J_4=J_5=0$ & 18 & \cellcolor{blue!25} $2 II^* + 4 I_1$ & \cellcolor{blue!25} $\{\mathbb{I}\}$ 
				& \cellcolor{blue!25} $H \oplus E_8(-1) \oplus E_8(-1)$ & \cellcolor{blue!25} $0$\\[0.1em]
				\hline
		\end{tabular}}
		\caption{Extensions of lattice polarizations for the \emph{standard} fibration}
		\label{fig:polarization_std}
	\end{subfigure} \\[0.5em]
	\begin{subfigure}{\textwidth}
		\centering
		\scalebox{\MyScale}{
			\begin{tabular}{| c || c  | c | c | c | c|}
				\hline 
				&&&&& \\[-0.9em]
				Fibration (\ref{eqn:alt_mod}) & $p_\mathcal{X}$ & Singular Fibers & $\operatorname{MW}(\pi_{\text{alt}}, \sigma_{\text{alt}})$ 
				& Lattice Polarization $\Lambda$ & $D(\Lambda)$ \\[0.1em]
				\hline\hline	
				&&&&& \\[-0.9em]
				generic & 16 & $I_8^* + 2 I_2 + 6 I_1$ & $\mathbb{Z}/2\mathbb{Z}$ 
				& $H \oplus E_7(-1) \oplus E_7(-1)$ & $\mathbb{Z}_2^2$\\[0.1em]
				\hline
				&&&&& \\[-0.9em]
				$\operatorname{Res}_t(D, \, E)=0$ & 16 & $I_8^* + III +  I_2 + 5 I_1$ & $\mathbb{Z}/2\mathbb{Z}$ 
				& $H \oplus E_7(-1) \oplus E_7(-1)$ & $\mathbb{Z}_2^2$\\[0.1em]
				\hline
				&&&&& \\[-0.9em]
				$\mathfrak{a}=0$ & 17 & $I_8^* + I_4 + 6 I_1$ & $\mathbb{Z}/2\mathbb{Z}$ 
				& $H \oplus E_8(-1) \oplus D_7(-1)$ & $\mathbb{Z}_4$\\[0.1em]
				\hline
				&&&&& \\[-0.9em]
				$J_{30}= 0$ & 17 & $I_8^* + 3 I_2 + 4 I_1$ & $\mathbb{Z}/2\mathbb{Z}$ 
				& $H \oplus E_7(-1) \oplus E_7(-1)  \oplus A_1(-1)$ & $\mathbb{Z}_2^3$\\[0.1em]
				\hline
				&&\cellcolor{red!25}&\cellcolor{red!25}&\cellcolor{red!25}&\cellcolor{red!25} \\[-0.9em]
				$J_4=0$ & 17 & \cellcolor{red!25} $I_{10}^* + I_2 + 6 I_1$ & \cellcolor{red!25} $\mathbb{Z}/2\mathbb{Z}$ 
				& \cellcolor{red!25} $H \oplus E_8(-1) \oplus E_7(-1)$ & \cellcolor{red!25} $\mathbb{Z}_2 $\\[0.1em]
				\hline
				&&\cellcolor{red!25}&\cellcolor{red!25}&\cellcolor{red!25}&\cellcolor{red!25} \\[-0.9em]
				$J_4=J_5=0$ & 18 & \cellcolor{red!25} $I_{12}^* + 6 I_1$ & \cellcolor{red!25} $\mathbb{Z}/2\mathbb{Z}$ 
				& \cellcolor{red!25} $H \oplus E_8(-1) \oplus E_8(-1)$ & \cellcolor{red!25} $0$\\[0.1em]
				\hline
		\end{tabular}}
		\caption{Extensions of lattice polarizations for the \emph{alternate} fibration}
		\label{fig:polarization_alt}
	\end{subfigure} \\[0.5em]
	\begin{subfigure}{\textwidth}
		\centering
		\scalebox{\MyScale}{
			\begin{tabular}{| c || c  | c | c | c | c|}
				\hline 
				&&&&& \\[-0.9em]
				Fibration (\ref{eqn:bfd_mod}) & $p_\mathcal{X}$ & Singular Fibers & $\operatorname{MW}(\pi_{\text{bfd}}, \sigma_{\text{bfd}})$ 
				& Lattice Polarization $\Lambda$ & $D(\Lambda)$ \\[0.1em]
				\hline\hline	
				&&&&& \\[-0.9em]
				generic & 16 & $II^* + I_2^* + 6 I_1$ & $\{\mathbb{I}\}$
				& $H \oplus E_8(-1) \oplus D_6(-1)$ & $\mathbb{Z}_2^2$\\[0.1em]
				\hline
				&&&&& \\[-0.9em]
				$\operatorname{Res}_t(t^{-2}F, \, t^{-3}G)=0$ & 16 & $II^* + I_2^* + II+ 4 I_1$ & $\{\mathbb{I}\}$ 
				& $H \oplus E_8(-1) \oplus D_6(-1)$ & $\mathbb{Z}_2^2$\\[0.1em]
				\hline
				&&&&& \\[-0.9em]
				$\mathfrak{a}=0$ & 17 & $II^* + I_3^* + 5 I_1$ & $\{\mathbb{I}\}$ 
				& $H \oplus E_8(-1) \oplus D_7(-1)$ & $\mathbb{Z}_4$\\[0.1em]
				\hline
				&&&&& \\[-0.9em]
				$J_{30} =0$ & 17 & $II^* + I_2^* + I_2 + 4 I_1$ & $\{\mathbb{I}\}$ 
				& $H \oplus E_8(-1) \oplus D_6(-1) \oplus A_1(-1)$ & $\mathbb{Z}_2^3$\\[0.1em]
				\hline
				&&\cellcolor{blue!25}&\cellcolor{blue!25}&\cellcolor{blue!25}&\cellcolor{blue!25} \\[-0.9em]
				$J_4=0$ & 17 & \cellcolor{blue!25}  $II^* + III^* + 5 I_1$ & \cellcolor{blue!25}  $\{\mathbb{I}\}$ 
				& \cellcolor{blue!25}  $H \oplus E_8(-1) \oplus E_7(-1)$ & \cellcolor{blue!25}  $\mathbb{Z}_2 $\\[0.1em]
				\hline
				&&\cellcolor{blue!25}&\cellcolor{blue!25}&\cellcolor{blue!25}&\cellcolor{blue!25} \\[-0.9em]
				$J_4=J_5=0$ & 18 & \cellcolor{blue!25} $2 II^* + 4 I_1$ & \cellcolor{blue!25} $\{\mathbb{I}\}$ 
				& \cellcolor{blue!25} $H \oplus E_8(-1) \oplus E_8(-1)$ & \cellcolor{blue!25} $0$\\[0.1em]
				\hline		
		\end{tabular}}
		\caption{Extensions of lattice polarizations for the \emph{base-fiber dual} fibration}
		\label{fig:polarization_bfd}
	\end{subfigure} \\[0.5em]
	\begin{subfigure}{\textwidth}
		\centering
		\scalebox{\MyScale}{
			\begin{tabular}{| c || c  | c | c | c | c|}
				\hline 
				&&&&& \\[-0.9em]
				Fibration (\ref{eqn:max_mod}) & $p_\mathcal{X}$ & Singular Fibers & $\operatorname{MW}(\pi_{\text{max}}, \sigma_{\text{max}})$ 
				& Lattice Polarization $\Lambda$ & $D(\Lambda)$ \\[0.1em]
				\hline\hline	
				&&&&& \\[-0.9em]
				generic & 16 & $I_{10}^* + 8 I_1$ & $\{\mathbb{I}\}$ 
				& $H \oplus D_{14}(-1)$ & $\mathbb{Z}_2^2$\\[0.1em]
				\hline
				&&&&& \\[-0.9em]
				$\mathrm{Res}_t(\alpha, \beta) = 0$, & 16 & $I_{10}^* + II + 6I_1$ & $\{\mathbb{I}\}$ & $H \oplus D_{14}(-1)$ & $\mathbb{Z}_2^2$\\[0.1em]
				\hline
				&&&&& \\[-0.9em]
				$\mathfrak{a}=0$ & 17 & $I_{11}^* + 7 I_1$ & $\{\mathbb{I}\}$ 
				& $H \oplus D_{15}(-1)$ & $\mathbb{Z}_4$\\[0.1em]
				\hline
				&&&&& \\[-0.9em]
				$J_{30}=0$ & 17 & $I_{10}^* + I_2 + 6 I_1$ & $\{\mathbb{I}\}$ 
				& $H \oplus D_{14}(-1) \oplus A_1(-1)$ & $\mathbb{Z}_2^3$\\[0.1em]
				\hline
				&&\cellcolor{red!25}&\cellcolor{red!25}&\cellcolor{red!25}&\cellcolor{red!25} \\[-0.9em]
				$J_4=0$ & 17 & \cellcolor{red!25} $I_{10}^* + I_2 + 6 I_1$ & \cellcolor{red!25} $\mathbb{Z}/2\mathbb{Z}$ 
				& \cellcolor{red!25} $H \oplus E_8(-1) \oplus E_7(-1)$ & \cellcolor{red!25} $\mathbb{Z}_2 $\\[0.1em]
				\hline
				&&\cellcolor{red!25}&\cellcolor{red!25}&\cellcolor{red!25}&\cellcolor{red!25} \\[-0.9em]
				$J_4=J_5=0$ & 18 & \cellcolor{red!25} $I_{12}^* + 6 I_1$ & \cellcolor{red!25} $\mathbb{Z}/2\mathbb{Z}$ 
				& \cellcolor{red!25} $H \oplus E_8(-1) \oplus E_8(-1)$ & \cellcolor{red!25} $0$\\[0.1em]
				\hline
		\end{tabular}}
		\caption{Extensions of lattice polarizations for the \emph{maximal} fibration}
		\label{fig:polarization_max}
	\end{subfigure} 
	\caption{Extensions of lattice polarization}
	\label{fig:polarization}
\end{figure}
\section{Classification of non-geometric heterotic models}
\label{sec:heterotic_models}
An eight-dimensional effective theory for the heterotic string compactified on $\mathbf{T}^2$ has a complex scalar field which takes its values in the Narain space \cite{MR834338}
\beq
\mathcal{D}_{2,18}/O(\Lambda^{2,18}),
\eeq
where $\mathcal{D}_{p, q}$ is the symmetric space for $O(p, q)$, i.e.,
\beq
\mathcal{D}_{p, q} = (O(p)\times O(q))\backslash O(p, q).
\eeq
The \emph{Narain space} is the quotient of the symmetric space for $O(2,18)$ by the automorphism group $O(\Lambda^{2,18})$ of the unique integral even unimodular lattice of signature $(2,18)$, i.e.,
\beq
\Lambda^{2,18}=H\oplus H \oplus E_8(-1)\oplus E_8(-1) \;.
\eeq
In an appropriate limit, the Narain space decomposes as a product of  spaces parameterizing the K\"ahler and complex structures on $\mathbf{T}^2$ as well as sixteen Wilson line expectation values around the two generators of $\pi_1(\mathbf{T}^2)$; see~\cite{MR867240} for details.  However, the decomposition is not preserved when the moduli vary arbitrarily. Families of heterotic models employing the full $O(\Lambda^{2,18})$ symmetry are therefore considered \emph{non-geometric} compactifications, because the K\"ahler and  complex structures on $\mathbf{T}^2$, and the Wilson line values, are not distinguished under the $O(\Lambda^{2,18})$-equivalences but instead are mingled together.
\par If we restrict ourselves to a certain index-two sub-group $O^+(\Lambda^{2,18}) \subset O(\Lambda^{2,18})$  in the construction above, the non-geometric models can be described by holomorphic modular forms. This is because the group $O^+(\Lambda^{2,18})$ is the maximal sub-group whose action preserves the complex structure on the symmetric space, and thus is the maximal sub-group for which modular forms are holomorphic.  The statement of the \emph{F-theory/heterotic string duality} in eight dimensions \cite{MR1403744} is the statement that the quotient space 
\beq
\label{k3moduli}
\mathcal{D}_{2,18}/ O^+(\Lambda^{2,18})
\eeq
coincides with the parameter space of elliptically fibered K3 surfaces with a section, i.e, the moduli space of F-theory models. This statement has been known in the mathematics literature as well; see, for example, \cite{MR2336040}. However, to construct the duality map between F-theory models and heterotic sting vacua explicitly,  one has to know the ring of modular forms relative to $O^+(\Lambda^{2,18})$ and their connection to the corresponding elliptically fibered K3 surfaces. However, this ring of modular forms is not known in general. We consider the restriction to a natural four-dimensional sub-space $\mathcal{D}_{2,4}/O^+(2,4;\mathbb{Z})$ of the space in Equation~(\ref{k3moduli}). Due to Theorem~\ref{thm:ROMF} the corresponding ring of modular forms is known.
\par Let $L^{2,4}$ be the lattice of signature $(2,4)$ which is the orthogonal complement of $E_7(-1)\oplus E_7(-1)$ in $\Lambda^{2,18}$. By insisting that the Wilson lines associated to the $E_7(-1)\oplus E_7(-1)$ sub-lattice are trivial, we restrict to heterotic vacua parameterized by the sub-space
\beq 
\mathcal{D}_{2,4}/O(L^{2,4}) \,.
\eeq
The corresponding degree-two cover is precisely the quotient space discussed above, namely
\beq \label{eqn:MSP+} \mathcal{D}_{2,4}/O^+(L^{2,4}) \,.
\eeq
For this natural four-dimensional sub-space in the full eighteen dimensional moduli space we will determine the duality map (and thus the quantum-exact effective interactions) between a dual F-theory/heterotic string pair in eight space-time dimensions. As we will show, the restriction to this sub-space describes the partial Higgsing of the corresponding heterotic gauge algebra $\mathfrak{g}=\mathfrak{e}_8 \oplus \mathfrak{e}_8$ or $\mathfrak{g}=\mathfrak{so}(32)$\ for the associated low energy supergravity theory. 
\par The Jacobian elliptic fibrations in Section~\ref{ssec:JAC-fibrations} provide the possible F-theories for $D=8$ compactifications over the four-dimensional sub-space in~(\ref{eqn:MSP+}) of the full eighteen dimensional moduli space. We now consider families of such non-geometric heterotic compactifications that naturally lead to compactifications for $D<8$. To start with, an inspection of our results from Section~\ref{ssec:JAC-fibrations} shows that for the dual F-theory models there are \emph{no} Jacobian elliptic fibrations on the sub-space~(\ref{eqn:MSP+}) with a Mordell-Weil group of positive rank. Non-torsion sections in a Weierstrass model are known to describe the charged matter fields of the corresponding F-theory model \cites{MR3083343, MR3270382}. Thus, we have the following:
\begin{corollary}
For generic families of non-geometric heterotic compactifications sampling the moduli space $\mathcal{D}_{2,4}/O^+(L^{2,4})$ there cannot be any charged matter fields.
\end{corollary}
\subsection{The \texorpdfstring{$\mathfrak{e}_8\oplus\mathfrak{e}_8$}{e8+e8}-string}
\label{ssec:e8+e8}
As we have seen in Section~\ref{ssec:modular}, the space in Equation~(\ref{eqn:MSP+}) parameterizes pseudo-ample K3 surfaces with $H \oplus E_7(-1)\oplus E_7(-1)$ lattice polarization. Section~\ref{sec:BFDFibration} shows that these K3 surfaces admit an elliptic fibration with section, one fiber of Kodaira type $I_2^*$ or worse, and another fiber of type precisely $II^*$. Here, we have used the lattice isomorphism
\beq
H \oplus E_7(-1)\oplus E_7(-1) \cong H \oplus E_8(-1)\oplus D_6(-1) \,.
\eeq 
Because of the presence of a $II^*$ fiber, the Mordell-Weil group is always trivial, including all cases with gauge symmetry enhancement. From a physics point of view as was argued in \cite{MR3366121}, assuming that one fiber is fixed and of Kodaira type $II^*$ will avoid ``pointlike instantons'' on the heterotic dual after further compactification to dimension six or below, at least for general moduli.
\par The key geometric fact for the construction of F-theory models is that Equation~(\ref{eqn:bfd_mod}) defines an elliptically fibered K3 surface $\mathcal{X}$ with section whose periods determine a point $\varpi \in \mathbf{H}_{2,2}$ up to the action of $\Gamma^+_{T}$, and with the coefficients in the defining equation being modular forms relative to $\Gamma^+_{T}$ of even characteristic. The explicit form of the F-theory/heterotic string duality on the moduli space in Equation~(\ref{eqn:MSP+}) then has two parts:  starting from $\varpi \in \mathbf{H}_{2,2}$ we always obtain a Jacobian elliptic fibration on the K3 surface $\mathcal{X}$ from Equation~(\ref{eqn:bfd_mod}).  Conversely, we can start with any Jacobian elliptic fibration given by the general equation
\beq
\label{eq:genform}
\begin{split}
	Y^2 & = X^3 + a \, t^2 X + b \, t^3 + c \, t^3 X + c \, d\, t^4  
	+ e\, t^4 X + (d \, e+ f) \, t^5 + g\, t^6 + t^7 \,.
\end{split}
\eeq
We then determine a point in $\varpi \in \mathbf{H}_{2,2}$ (up to the action of $\Gamma^+_{T}$) by calculating the periods of the holomorphic two-form $\omega_{\mathcal{X}}=dt\wedge dX/Y$ over a basis of the lattice $H \oplus E_7(-1)\oplus E_7(-1)$ in $H^2(\mathcal{X},\mathbb{Z})$. It follows that for some non-vanishing scale factor $\lambda$ we have
\beq
\label{params}
\begin{split}
	c= - \lambda^{10} J_5(\varpi),\;
	d= - \frac{\lambda^8}{3}J_4(\varpi),\;
	e = -3 \, \lambda^4 J_2(\varpi), \;
	f =  \lambda^{12} J_6(\varpi),\;
	g= -2 \, \lambda^6 J_3(\varpi) , 
\end{split}
\eeq
and $a = - 3\, d^2, \quad b = -2 \, d^3$. Under the restriction of $\mathbf{H}_{2,2}$ to the Siegel upper half-plane $\mathbb{H}_2$, we have $d=0$ and
\beq
\label{SiegelRestriction2}
[J_2(\varpi) : J_3(\varpi) : J_4(\varpi) : J_5(\varpi) : J_6(\varpi)] =
[\psi_4(\tau) : \psi_6(\tau) : 0 : 2^{12} 3^5 \chi_{10}(\tau) : 2^{12} 3^6 \chi_{12}(\tau)] \,,
\eeq
as points in the four dimensional weighted projective space $\mathbb{W}\mathbb{P}(2,3,4,5,6)$, where $\psi_4, \psi_6, \chi_{10},$ and $\chi_{12}$ are Siegel modular forms of respective weights $4, 6, 10$ and $12$ introduced by Igusa in \cite{MR0141643}.  Moreover, a simple rescaling reduces Equation~(\ref{eqn:bfd_mod}) to
\beq
\label{eqn:std_red}
Y^2 = X^3 - t^3  \Big(\frac{1}{48} \psi_4(\tau) t + 4 \chi_{10}(\tau) \Big) X + t^5 \Big(t^2 - \frac{1}{864} \psi_6(\tau) \, t + \chi_{12}(\tau) \Big) \,,
\eeq
which is precisely the equation derived in \cite{MR3366121} for the F-theory dual of a non-geometric heterotic theory with gauge algebra $\mathfrak{g} = \mathfrak{e}_8 \oplus \mathfrak{e}_7$ and one non-vanishing Wilson line parameter.  We have proved the following:
\begin{proposition}
\label{prop:het_vacua1a}
Equation~(\ref{eq:genform}) defines the F-theory dual of a non-geometric heterotic theory with gauge algebra $\mathfrak{g} = \mathfrak{e}_8 \oplus \mathfrak{so}(12)$.
\end{proposition}
\subsection{Condition for five-branes and supersymmetry}
\label{sec:susy}
The strategy for constructing \emph{families} of non-geometric heterotic compactifications is the following:  start with a compact manifold $\mathfrak{Z}$ as parameter space and a line bundle  $\Lambda  \to \mathfrak{Z}$.  Choose sections $c(z)$, $d(z)$, $e(z)$, $f(z)$, and $g(z)$ of the bundles $\Lambda^{\otimes 10}$, $\Lambda^{\otimes 8}$,  $\Lambda^{\otimes 4}$, $\Lambda^{\otimes 12}$, and $\Lambda^{\otimes 6}$, respectively; then, for each point $z\in \mathfrak{Z}$, there is a non-geometric heterotic compactification given by Equation~(\ref{eq:genform}) with $c=c(z)$, $d=d(z)$, etc., and $a = - 3\, d(z)^2$, $b = -2 \, d(z)^3$ and moduli $\varpi \in \mathbf{H}_{2,2}$ and $O^+(L^{2,4})$ symmetry such that Equations~(\ref{params}) hold. 
\par Appropriate five-branes must still be inserted on $\mathfrak{Z}$ as dictated by the geometry of the corresponding family of K3 surfaces. The change in the singularities and the lattice polarization for the fibration~(\ref{eqn:bfd_mod}) occur along three loci of co-dimension one, namely, $\mathfrak{a}=0$, $J_{30}=0$, and $J_4=0$; see Section~\ref{sec:confluences}. Each locus is the fixed locus of elements in $\Gamma_{T} \, \backslash \, \Gamma^{+}_{T}$. It is trivial to write down the reflections in $O^+(L^{2,4}) \backslash SO^+(L^{2,4})$ corresponding to $\mathfrak{a}=0$, $J_{30}=0$, and $J_4=0$, respectively; see \cite{MR4015343}.
\par From the point of view of K3 geometry, given as a reflection in a lattice element $\delta$ of square $-2$ we have the following: if the periods are preserved by the reflection in $\delta$, then $\delta$ must belong to the N\'eron-Severi lattice of the K3 surface.  That is, the N\'eron-Severi lattice is enlarged by adjoining $\delta$. We already showed in Section~\ref{sec:confluences} that there are three ways an enlargement can happen: the lattice $H \oplus E_7(-1) \oplus E_7(-1)$ of rank sixteen can be extended to $H\oplus E_7(-1)\oplus E_7(-1) \oplus \langle -2 \rangle$, $H \oplus E_8(-1)\oplus E_7(-1)$, or $H \oplus E_8(-1) \oplus D_7(-1)$, each of rank seventeen. 
\par On the heterotic side these five-brane solitons are easy to see -- we derived the corresponding confluences of singular fibers in Section~\ref{sec:confluences}:  when $J_{30}=0$, we have a gauge symmetry enhancement from $\mathfrak{e}_8 \oplus \mathfrak{so}(12)$ to include an additional $\mathfrak{su}(2)$, and the  parameters of the theory include a Coulomb branch for that gauge theory on which the Weyl group $W_{\mathfrak{su}(2)}=\mathbb{Z}_2$ acts.  Thus, there is a five-brane solution in which the field has a $\mathbb{Z}_2$ ambiguity encircling the location in the moduli space of enhanced gauge symmetry.  When $J_4 =0$, we have an enhancement to $\mathfrak{e}_8 \oplus \mathfrak{e}_7$ gauge symmetry, and, when $\mathfrak{a}=0$ an enhancement to $\mathfrak{e}_8 \oplus \mathfrak{so}(14)$.  Further enhancement to $\mathfrak{e}_8 \oplus \mathfrak{e}_8$ gauge symmetry occurs along $J_4=J_5=0$.
\par To understand when such families of compactifications are supersymmetric, we adapt the discussion in \cite{MR3366121}: a heterotic compactification on $\mathbf{T}^2$ with parameters given by $\varpi \in \mathbf{H}_{2,2}$ is dual to the F-theory compactification on the elliptically fibered K3 surface $\mathcal{X}(\varpi)$.  For sections $c(z)$, $d(z)$, $e(z)$, $f(z)$, and $g(z)$ of line bundles over $\mathfrak{Z}$, we have a criterion for when F-theory compactified on the elliptically fibered manifold~(\ref{eq:genform}) is supersymmetric: this is the case if and only if the total space defined by Equation~(\ref{eq:genform}) -- now considered as an elliptic fibration over a base space locally given by variables $t$ and $z$ -- is itself a Calabi--Yau manifold.  The base space of the elliptic fibration is a $\mathbb{P}^1$-bundle $\pi:\mathfrak{W} \to \mathfrak{Z}$ which takes the form $\mathfrak{W}=\mathbb{P}(\mathcal{O}\oplus \mathcal{M})$ where $\mathcal{M}\to \mathfrak{Z}$ is the normal bundle of $\Sigma_0:=\{t=0\}$ in $\mathfrak{W}$. Monomials of the form $t^n$ are then considered sections of the line bundles $\mathcal{M}^{\otimes n}$. We also set $\Sigma_\infty:=\{t=\infty\}$ such that $-K_{\mathfrak{W}} = \Sigma_0 + \Sigma_\infty + \pi^{-1}(-K_{\mathfrak{Z}})$.
\par When the elliptic fibration~(\ref{eq:genform}) is written in Weierstrass form, the coefficients of $X^1$ and $X^0$ must again be sections of $\mathcal{L}^{\otimes 4}$ and $\mathcal{L}^{\otimes 6}$, respectively, for a line bundle $\mathcal{L} \to \mathfrak{W}$. The condition for supersymmetry of the total space is $\mathcal{L}=\mathcal{O}_{\mathfrak{W}}(-K_{\mathfrak{W}})$. Restricting the various terms in Equation~(\ref{eq:genform}) to $\Sigma_0$, we find relations
\beq
\label{eqn:bundle_relations}
\begin{split}
	(\mathcal{L}|_{\Sigma_0})^{\otimes 4} & = \Lambda^{\otimes 4} \otimes \mathcal{M}^{\otimes 4}  = \Lambda^{\otimes 10} \otimes \mathcal{M}^{\otimes 3}
	=\Lambda^{\otimes 16} \otimes \mathcal{M}^{\otimes 2},\\
	(\mathcal{L}|_{\Sigma_0})^{\otimes 6} & = \mathcal{M}^{\otimes 7} =\Lambda^{\otimes 6} \otimes \mathcal{M}^{\otimes 6}\\
	&=\Lambda^{\otimes 12} \otimes \mathcal{M}^{\otimes 5} =\Lambda^{\otimes 18} \otimes \mathcal{M}^{\otimes 4} =\Lambda^{\otimes 24} \otimes \mathcal{M}^{\otimes 3} .
\end{split}
\eeq
Thus, it follows that  $\mathcal{M}=\Lambda^{\otimes 6}$ and $\mathcal{L}|_{\Sigma_0}=\Lambda^{\otimes 7}$ (up to torsion) and the $\mathbb{P}^1$-bundle takes the form $\mathfrak{W}=\mathbb{P}(\mathcal{O}\oplus \Lambda^{\otimes 6})$. Since $\Sigma_0$ and $\Sigma_\infty$ are disjoint, the condition for supersymmetry is equivalent to $\Lambda = \mathcal{O}_{\mathfrak{Z}}(-K_{\mathfrak{Z}})$. We have proved the following:
\begin{proposition}
\label{prop:het_vacua1b}
Equation~(\ref{eq:genform}) defines a supersymmetric family of non-geometric heterotic vacua with gauge algebra $\mathfrak{e}_8 \oplus \mathfrak{so}(12)$  over  a  compact parameter space $\mathfrak{Z}$ equipped with the line bundle  $\Lambda = \mathcal{O}_{\mathfrak{Z}}(-K_{\mathfrak{Z}}) \to \mathfrak{Z}$ if $c(z)$, $d(z)$, $e(z)$, $f(z)$, and $g(z)$ in Equation~(\ref{params}) are sections of the bundles $\Lambda^{\otimes 10}$, $\Lambda^{\otimes 8}$,  $\Lambda^{\otimes 4}$, $\Lambda^{\otimes 12}$, and $\Lambda^{\otimes 6}$, respectively.
\end{proposition}
\subsection{Double covers and pointlike instantons}
To a reader familiar with elliptic fibrations, it might come as a surprise that the Weierstrass model we considered in Equation~(\ref{eq:genform}) did not simply have two fibers of Kodaira type $III^*$ and a trivial Mordell-Weil group. On each K3 surface endowed with a $H \oplus E_7(-1) \oplus E_7(-1)$ lattice polarization such a fibration exists, and we constructed it in Equation~(\ref{eqn:std_mod}). However,  it is not guaranteed that the fibration extends across any parameter space, and there might be anomalies present. 
\par Starting from $\varpi \in \mathbf{H}_{2,2}$ we always obtain a pair of Jacobian elliptic fibrations on the K3 surface $\mathcal{X}$ from Equation~(\ref{eqn:std_mod}).  The indeterminacy of the sign $\pm \mathfrak{a}$ forces us to consider a pair of fibrations related by a holomorphic, symplectic isomorphism. Conversely, we can start with the general Jacobian elliptic fibration, normalized with $\delta=1$ and given by the general equation
\beq
\label{eq:genform1b}
\begin{split}
	Y^2 & = X^3 -\varepsilon t^3 X -3 \alpha t^4 X + \zeta t^5  -\gamma t^5 X  -2 \beta t^6 +  t^7  \,.
\end{split}
\eeq
If we then determine a point in $\varpi \in \mathbf{H}_{2,2}$ by calculating the periods of the holomorphic two-form $\omega_{\mathcal{X}}=dt\wedge dX/Y$ over a basis of the lattice $H \oplus E_7(-1)\oplus E_7(-1)$ in $H^2(\mathcal{X},\mathbb{Z})$, we must have $J_6\not = 0$ since otherwise the terms $t^7$ and $\zeta t^5$ could not both be present. It follows that for some non-vanishing scale factor $\lambda$ we have
\beq
\label{params2}
\begin{split}
\alpha = \lambda^4 J_2(\varpi) \,, \quad
\beta =  \lambda^6 J_3(\varpi) \,, \quad
\zeta = \lambda^{12}J_6(\varpi) \,,\\
\epsilon = \lambda^{10} \frac{(J_5 \mp \mathfrak{a})(\varpi)}{2} \,, \quad
\gamma = \lambda^{-2}   \frac{(J_5 \pm \mathfrak{a})(\varpi)}{2 J_6(\varpi)}\,, \quad
\end{split}
\eeq
such that  $\gamma \epsilon =  \lambda^{8}J_4(\varpi)$ and $\gamma \zeta + \epsilon = \lambda^{10}J_5(\varpi)$.
\par In order to construct \emph{families} of non-geometric compactifications, we vary the heterotic vacua over a parameter space $\mathfrak{Z}$ as in Section~\ref{sec:susy}, the functions $\alpha, \beta, \gamma\varepsilon, \zeta$ are again sections of line bundles  $\Lambda^{\otimes 2k}  \to \mathfrak{Z}$ for $k=2, 3, 4, 6$.  The condition for supersymmetry already established in Section~\ref{sec:susy} yields $\Lambda = \mathcal{O}_{\mathfrak{Z}}(-K_{\mathfrak{Z}})$. We want to take the square root of the line bundle $\Lambda^{\otimes 20}=\Lambda^{\otimes 8} \otimes  \Lambda^{\otimes 12}$, that is, construct a line bundle $\mathcal{N} \to \mathfrak{Z}$ with $\mathcal{N}^{\otimes 2} = \Lambda^{\otimes 20}$ such that $\mathfrak{a}$ becomes a section of the new line bundle $\mathcal{N}$.
\par If the line bundle $\mathcal{N}$ is effective, i.e., $\mathcal{N} \cong \mathcal{O}_{\mathfrak{Z}}(\mathrm{D})$ for some smooth, effective divisor $\mathrm{D}$ in $\mathfrak{Z}$ -- which is satisfied if $\dim{H^0(\mathfrak{Z}, \mathcal{N})} >0$ --  then a double cover $\varphi: \mathfrak{Y} \to \mathfrak{Z}$ can be constructed that is ramified over $\mathrm{D}$. The double cover $\varphi: \mathfrak{Y} \to \mathfrak{Z}$ is defined in terms of the line bundle $\mathcal{N} \to \mathfrak{Z}$ and a non-trivial section $\sigma$ of $\mathcal{N}^{\otimes 2}$ as follows:  (1) $\mathfrak{Z}$ is embedded in the total space of $\mathcal{N}^{\otimes 2}$ and given locally by an equation of the form $z = \sigma$, (2)  $\mathfrak{Y}$ is embedded in the total space of $\mathcal{N}$ and given locally by the equation $y^2 = \sigma$, (3) the double cover $\varphi$ is the restriction of the square map $\mathcal{N} \to \mathcal{N}^{\otimes 2}$, (4) the branch locus is given by $\{ \sigma = 0 \} \subset \mathfrak{Z}$. In turn, it is known that if $H_1(\mathfrak{Z})=0$ then the double cover $\varphi: \mathfrak{Y} \to \mathfrak{Z}$ is uniquely determined by its branch locus; see \cite{MR2030225}. In our situation, we set $\sigma =\mathfrak{a}^2 = J_5^2 - 4 J_4 J_6$.
\par Thus, Equation~(\ref{eq:genform1b}) can define a supersymmetric family of non-geometric heterotic vacua with gauge algebra $\mathfrak{e}_7 \oplus \mathfrak{e}_7$ over a compact parameter space $\mathfrak{Y}$, if $\mathfrak{Y}$ is obtained as the double cover of $\mathfrak{Z}$ in Section~\ref{sec:susy}  branched along $\mathfrak{a}=0$. We showed in Section~\ref{sec:confluences} that for $J_4=0$ we have $J_5^2=\mathfrak{a}^2$ and the choice of square root $\mathfrak{a}=\pm J_5$ determines which of the two fibers of type $III^*$ is extended to a fiber of type $II^*$; Equation~(\ref{eqn:std_mod}) then reduces to Equation~(\ref{eqn:std_red}) which is precisely the equation derived in \cite{MR3366121} for the F-theory dual of a heterotic theory with gauge algebra $\mathfrak{g} = \mathfrak{e}_8 \oplus \mathfrak{e}_7$ and one non-vanishing Wilson line parameter. 
\par However, the underlying K3 surface in Equation~(\ref{eqn:std_mod}) acquires a singularity whenever $J_6 \to 0$.  For the coefficient of $t^5 X$ in Equation~(\ref{eqn:std_mod}) to remain well defined, we may require $J_6 \not =0$ over $\mathfrak{Z}$ which implies that $J_6$ is a trivializing section for the bundle $\Lambda^{\otimes 12}$; in particular, we have $\Lambda^{\otimes 12} \cong \mathcal{O}_{\mathfrak{Z}}$.  Then, the conditions derived for supersymmetry and the existence of the double cover are similar to the conditions governing global and local anomaly cancellation in \cites{MR2854198, MR2815730}.  On the other hand, if a general family of non-geometric compactifications in Equation~(\ref{eq:genform1b}) intersects the locus $J_6 = 0$, the K3 surface in Equation~(\ref{eqn:std_mod}) becomes singular while at the same time the curvature  of the bundle $\Lambda$ also gets concentrated in an infinitesimal region of the base space because of Equation~(\ref{params2}). In \cites{MR1457988} the authors identified these points as \emph{point-like instantons} in the heterotic string, which correspond to the curvature acquiring singularities at singular points on the (blow-down of the) K3 surface. 
\subsection{The \texorpdfstring{$\mathfrak{so}(32)$}{so(32)}-string}
\label{sec:so32}
In Section~\ref{sec:family_of_K3s} we showed that a K3 surface $\mathcal{X}$ with lattice polarization $H \oplus E_7(-1) \oplus E_7(-1)$ also admits two additional  fibrations, which we called the  \emph{alternate} and the \emph{maximal} Jacobian elliptic fibration. These turn out to be related to the $\mathfrak{so}(32)$ heterotic string. 
\par We will now establish the explicit form of the F-theory/heterotic string duality on the moduli space~(\ref{eqn:MSP+}) for this case. The intrinsic property of the elliptically fibered K3 surfaces which lead to the corresponding F-theory backgrounds is the requirement that there is one singular fiber in the fibration of type $I_n^*$ for some $n\ge 8$ and a two-torsion element in the Mordell-Weil group. Then, a slight modification of the argument in \cite{MR1457988}*{Sec.~4} or \cite{MR2826187}*{App.~A.} shows that the corresponding Weierstrass equation takes the form
\beq
\label{eq:genform2}
Y^2  = X^3  + \Big(t^3 +e \, t + g\Big) \, X^2 + \Big(-3 \, d \, t^2+c \, t+f\Big) \, X \,.
\eeq
In the language of string theory, we obtain an F-theory model on a K3 surface, elliptically fibered over $\mathbb{P}^1$, with non-zero flux of an antisymmetric two-form $B$ through the sphere.  Such F-theory compactifications were first analyzed by Witten in \cite{MR1615617} in the limiting locus when the elliptic fibration becomes isotrivial and also discussed in \cites{MR1658283, MR1614531}. The picture was later extended to general F-theory elliptic fibrations in \cite{MR1797021}. Only the cohomology class of the antisymmetric two-form $B$ has a physical meaning. Accordingly, the value of the flux is quantized  and it is fixed to be equal to $\omega/2$ where $\omega$ denotes the K\"ahler class of the sphere.  In terms of the Jacobian elliptic fibration, the non-zero flux is generated by the non-trivial two-torsion element $(X,Y)=(0,0)$ of the Mordell-Weil group.
\par The explicit form of the F-theory/heterotic string duality on the moduli space in Equation~(\ref{eqn:MSP+}) for the $\mathfrak{so}(32)$ heterotic string then has two parts:  starting from a period point $\varpi \in \mathbf{H}_{2,2}$ and the action of $\Gamma^+_{T}$, we always obtain a Jacobian elliptic fibration with non-trivial two-torsion element in the Mordell--Weil group on the K3 surface $\mathcal{X}$  from Equation~(\ref{eqn:bfd_mod}).  Conversely, we can start with any Jacobian elliptic fibration in Equation~(\ref{eq:genform2}). We then determine a point in $\mathcal{D}_{2,4}$ by calculating the periods of the holomorphic two-form $\omega_{\mathcal{X}}$ over a basis of the period lattice  $H \oplus E_7(-1)\oplus E_7(-1)$ in $H^2(\mathcal{X},\mathbb{Z})$, such that a change of marking corresponds to the action of a modular transformation in $\Gamma^+_{T}$.  As before, it follows that for some non-vanishing scale factor $\lambda$ Equations~(\ref{params}) must hold. The gauge algebra is enhanced to $\mathfrak{so}(24) \oplus \mathfrak{su}(2)^{\oplus 2}$. It follows as in \cites{MR1416960,MR1643100} that the gauge group of this model is $(\operatorname{Spin}(24)\times SU(2)\times SU(2))/\mathbb{Z}_2$.  Thus, we have proved:
\begin{proposition}
\label{prop:het_vacua2a}
Equation~(\ref{eq:genform2}) defines the F-theory on an elliptic K3 surface with non-zero flux of an antisymmetric two-form $B$ through the sphere that is dual to a non-geometric heterotic theory with gauge algebra $\mathfrak{g} = \mathfrak{so}(24) \oplus \mathfrak{su}(2)^{\oplus 2}$.
\end{proposition}
The condition for five-branes and supersymmetry for families of non-geometric heterotic compactifications is completely analogous to Section~\ref{sec:susy}, and their construction is easily carried out as Equation~(\ref{eqn:alt_mod}) establishes the connection between the parameters and the modular forms relative to $\Gamma^+_{T}$. In particular, for $J_{30}=0$, we find a gauge symmetry enhancement to include an additional $\mathfrak{su}(2)$; when $\mathfrak{a}=0$ the gauge algebra is enhanced to $\mathfrak{so}(24) \oplus \mathfrak{su}(4)$. This follows from the results in Section~\ref{sec:confluences}.
\medskip

\par  In the mathematical classification of all distinct Jacobian elliptic fibrations supported on the family of K3 surfaces $\mathcal{X}$ in Equation~(\ref{quartic1}) for Picard rank 16 there is a fundamental difference compared to the classification in Picard rank $17$ and $18$. Here, the family of K3 surfaces $\mathcal{X}$ in Equation~(\ref{quartic1}) admits the additional Jacobian elliptic fibration given in Equation~(\ref{eqn:max_mod}) which we called the maximal fibration. In fact, the general elliptic fibration that has a singular fiber of Kodaira type $I_{8}^*$ over $t=\infty$ is of the form
\beq
\label{eq:genform2b}
Y^2 = X^3 + \Big( a' t^3 + b' t^2 + c' t + d'\Big) X^2 + \Big(e' t^2 + f' t + g'\Big) X + \Big( h' t + k' \Big) \,.
\eeq
The F-theory model determined by Equation~(\ref{eq:genform2b}) does \emph{not} admit a two-torsion element in the Mordell-Weil group and is polarized by the lattice $H \oplus E_7(-1)\oplus E_7(-1)$.  However, for $J_4=0$,  the two cases in Equation~(\ref{eq:genform2}) and~(\ref{eq:genform2b}) coincide. That is, after using Equation~(\ref{SiegelRestriction}) and a simple rescaling, both Equation~(\ref{eqn:alt_mod}) and Equation~(\ref{eqn:max_mod}) restrict to
\beq
\label{eqn:alt_red}
Y^2  = X^3  + \Big(t^3- \frac{1}{48} \, \psi_4(\tau) \, t- \frac{1}{864} \, \psi_6(\tau) \Big)  X^2 - \Big(4 \, \chi_{10}(\tau) \, t-\chi_{12}(\tau)\Big) X \,,
\eeq
which is precisely the equation derived in \cite{MR3366121} for the F-theory dual of a heterotic theory with gauge algebra $\mathfrak{g} = \mathfrak{so}(28) \oplus \mathfrak{su}(2)$ and one non-vanishing Wilson line parameter.  In  the limit $J_4 \to 0$, two nodes in Equation~(\ref{eq:genform2b}) coalesce and form a fiber of Kodaira type $I_2$ that generates the $ \mathfrak{su}(2)$-gauge enhancement. At the same time, when the coefficients $h', k' \to 0$ vanish, a non-trivial two-torsion element is generated in the Mordell-Weil group. We have proved the following:
\begin{proposition}
Equation~(\ref{eqn:max_mod}) defines the F-theory on an elliptic K3 surface dual to a non-geometric heterotic theory with gauge algebra $\mathfrak{g} = \mathfrak{so}(28)$.
\end{proposition}
The condition for five-branes and supersymmetry for families of non-geometric heterotic compactifications in this second case is again completely analogous to Section~\ref{sec:susy}, and their construction is easily carried out since Equation~(\ref{eqn:max_mod}) establishes the connection between the parameters and the modular forms relative to $\Gamma^+_{T}$.  In particular, it was shown in Section~\ref{sec:confluences} that for $\mathfrak{a}=0$, we have a gauge symmetry enhancement to $\mathfrak{so}(30)$, and for $J_4 = J_5 = 0$ the gauge group is further enhanced to $\operatorname{Spin}(32)/\mathbb{Z}_2$.
\section{Discussion and Outlook}
We classified all non-geometric heterotic models obtained by the partial Higgsing using two Wilson lines of the heterotic gauge algebra $\mathfrak{g}=\mathfrak{e}_8 \oplus \mathfrak{e}_8$ or $\mathfrak{so}(32)$ for the associated low energy effective eight-dimensional supergravity theory. The surprising result: as opposed to the Higgsing of the heterotic gauge algebra using only one Wilson lines, the Higgsing with two Wilson lines produces two different branches for each type of heterotic string theory. We interpreted the dual F-theory models as Jacobian elliptic fibrations supported on the family of K3 surfaces with canonical $H \oplus E_7 (-1) \oplus E_7 (-1)$ polarization. The inequivalent Jacobian elliptic fibrations were classified, and we found the defining equations for Weierstrass models whose parameters are modular forms generalizing well-known Siegel modular forms. Two of these fibrations correspond to the Higgsing of heterotic gauge algebra $\mathfrak{g}=\mathfrak{e}_8 \oplus \mathfrak{e}_8$ to either $\mathfrak{e}_8 \oplus \mathfrak{so}(12)$ or $\mathfrak{e}_7 \oplus \mathfrak{e}_7$, the other two are related to the Higgsing of $\mathfrak{g}=\mathfrak{so}(32)$ to either $\mathfrak{so}(24) \oplus \mathfrak{su}(2)^{\oplus 2}$ or $\mathfrak{so}(28)$. In the former case, the two fibrations are differentiated by avoiding or supporting  ``pointlike instantons'' on the heterotic dual after further compactification to dimension six. In the latter case, the two fibrations are distinguished by trivial or non-trivial flux of an antisymmetric two-form $B$ through the base of the elliptic fibration.  We demonstrated how the fibrations can be used to construct \emph{families} of non-geometric heterotic compactifications. The necessary condition for five-branes and supersymmetry were determined explicitly. Therefore, our results provide a significant generalization of existing results in \cites{MR3417046, MR2826187, MR3366121, MR3712162}.
\par As a result, this article provides a complete description of the F-theory/heterotic string duality in $D=8$ with two Wilson lines. For no or one non-trivial Wilson line parameters, an analogous approach has been proven to also provide a quantum-exact effective description of non-geometric heterotic models \cites{MR2826187,MR3366121,MR3417046}. We expect that the analysis carries over in our case. Since there is no microscopic description of the dual F-theory, the explicit form of the F-theory/heterotic string duality in this article also provides new insights into the physics of F-theory compactifications.  One of the conclusion of the work in \cite{MR3366121} was that taking a heterotic compactification even a ``small distance'' from the large radius limit destroys the traditional semiclassical interpretation and no longer allows us to discuss the compactification as being that of a manifold with a bundle. We point out that this is not unlike what happens in type II compactifications, where the analysis of $\Pi$-stability \cites{MR1915365, MR1840318} shows that going any distance away from large radius limit, no matter how  small, necessarily changes the stability conditions on some D-brane classes and so destroys the semiclassical interpretation of the theory. Thus, we expect that our results might be of importance for a better understanding of non-perturbative aspects of the heterotic string, for example, as it relates to NS5-branes states and small instantons \cites{MR1381611,MR1769478}. 
\par Moreover, we expect that the close connection between modular forms and equations presented in this article will enable us to use a well known F-theory construction and build interesting classes of non-geometric heterotic compactifications which have duals described in terms of K3-fibered Calabi-Yau manifolds. The starting point is the heterotic string compactified on a torus, and utilizes the non-perturbative duality symmetries which this theory possesses. The construction was explained in considerable detail in \cites{MR2826187, MR3366121} and was then used to obtain many examples realizing families of non-geometric heterotic compactifications  with one Wilson line parameter \cites{Mayrhofer_2017, MR3562156,MR3933163} using the classification of pencils of genus-two curves by Namikawa and Ueno \cites{MR0384794, MR0369362, MR0319996}. In \cite{Kimura2019UnbrokenE} our classification result was already used to construct certain examples of families of non-geometric heterotic string vacua with two Wilson line parameters and corresponding Calabi-Yau threefolds. It is interesting to ask whether a systematic program can be carried out, classifying all resulting compactifications in six dimensions. We leave this question for future work.
\bibliographystyle{amsplain}
\bibliography{ref}{}
\end{document}